\documentclass[a4paper]{article}
\usepackage{amsmath}
\newcommand{\be}{\begin{equation}}
\newcommand{\ee}{\end{equation}}
\newcommand{\bea}{\begin{eqnarray*}}
\newcommand{\eea}{\end{eqnarray*}}

\newcommand{\D}{{\cal D}}
\newcommand{\A}{{\cal A}}
\newcommand{\Li}{{\cal L}}
\newcommand{\s}{\sigma}
\newcommand{\p}{\partial}

\newcommand{\la}{\langle}
\newcommand{\ra}{\rangle}
\newcommand{\raa}{\rightarrow}
\newcommand{\ci}{C^{\infty}}
\newcommand{\C}{{\cal C}}
\newcommand{\E}{\ell}
\newcommand{\m}{\!\!\mid}
\newcommand{\cl}{{\cal L}}
\newcommand{\N}{{\mathbf N}}
\newcommand{\Z}{{\mathbf Z}}
\newcommand{\R}{{\mathbf R}}
\newcommand{\F}{{\cal F}}

\begin{document}

\title{Automorphisms of\\ quantum and classical Poisson algebras
\footnote{This work was partially supported by NATO Grant
RD/DOC/020204-3 and by MCESR Grant RD/C.U.L./02-010}}
\author{J. Grabowski, N. Poncin}\maketitle

\newtheorem{re}{Remark}
\newtheorem{theo}{Theorem}
\newtheorem{prop}{Proposition}
\newtheorem{lem}{Lemma}
\newtheorem{cor}{Corollary}
\newtheorem{ex}{Example}

\begin{abstract}
We prove Pursell-Shanks type results for the Lie algebra $\D(M)$
of all linear differential operators of a smooth manifold $M$, for
its Lie subalgebra $\D^1(M)$ of all linear first-order
differential operators of $M$, and for the Poisson algebra
$S(M)=Pol(T^*M)$ of all polynomial functions on $T^*M,$ the
symbols of the operators in $\D(M).$

Chiefly however we provide explicit formulas describing completely
the automorphisms of the Lie algebras $\D^1(M),$ $S(M),$ and
$\D(M).$
\end{abstract}

\section{Introduction}

The classical result of Pursell and Shanks \cite{PS}, which states
that the Lie algebra of smooth vector fields of a smooth manifold
characterizes the smooth structure of the variety, is the starting
point of a multitude of papers.

There are similar results in particular geometric situations---for
instance for hamiltonian, contact or group invariant vector
fields---for which specific tools have each time been constructed,
\cite{O,A,AG,HM}, in the case of Lie algebras of vector fields
that are modules over the corresponding rings of functions,
\cite{Am,JG,S}, as well as for the Lie algebra of (not leaf but)
foliation preserving vector fields, \cite{JG3}.

The initial objective of the present paper was to prove that the
Lie algebra $\D(M)$ of all linear differential operators
$D:\ci(M)\raa\ci(M)$ of a smooth manifold $M$, determines the
smooth structure of $M$. Beyond this conclusion, we present a
description of all automorphisms of the Lie algebra $\D(M)$ and
even of the Lie subalgebra $\D^1(M)$ of all linear first-order
differential operators of $M$ and of the Poisson algebra
$S(M)=Pol(T^*M)$ of polynomial functions on the cotangent bundle
$T^*M$ (the symbols of the operators in $\D(M)$), the
automorphisms of the two last algebras being of course canonically
related with those of $\D(M)$. In each situation we obtain an
explicit formula, for instance---in the case of $\D(M)$---in terms
of the automorphism of $\D(M)$ implemented by a diffeomorphism of
$M$, the conjugation-automorphism of $\D(M)$, and the automorphism
of $\D(M)$ generated by the derivation of $\D(M)$ associated to a
closed $1$-form of $M$.

In the first part of our work, the approach is purely algebraic.
In Section $2$, we heave $\D(M)$ and $S(M)$ on a general algebraic
level and define the notions "quantum Poisson algebra" $\D$ resp.
"classical Poisson algebra" $S$, classical limit of $\D$. In
Section $3,$ we show that if two (quantum or classical) Poisson
algebras are isomorphic as Lie algebras, their "basic algebras of
functions" are isomorphic as associative algebras---an algebraic
Shanks-Pursell type result, which naturally implies our previously
described initial goal. The leading idea of the proof is the
algebraic characterization, under a minimal condition, of
functions as those $D\in\D$ or $P\in S$ for which $ad_D$, resp.
$ad_P$, is locally nilpotent.

In the second part of the article, we switch to the concrete
geometric context. In this introduction, we will confine ourself
to a very rough description of the quite technical computations of
Section $7$ that give all automorphisms of $\D(M),$ calculations
based upon the result in the $S(M)$-case, Section $6$, itself
founded on the $\D^1(M)$-case, Section $5$. The utilization---in
addition to the just mentioned algebra hierarchy---of the
preliminary detected conjugation- and derivation-automorphisms,
Section $4,$ and the suitable use of the normal ordering method
(i.e. the local polynomial representation of differential
operators), allow to reduce the problem to the determination of
intertwining operators between some modules of the Lie algebra of
vector fields and to conclude.

\section{Definitions and tools}

By a \textit{quantum Poisson algebra} we understand an associative
filtered algebra ${\cal D}={\cup}_{i=0}^{\infty}{\cal D}^{i}$,
$\D^{i}\subset\D^{i+1}$, $\D^{i}\cdot\D^{j}\subset\D^{i+j}$ (where
$\cdot$ denotes the multiplication of ${\cal D}$), with unit $1$
over a field $K$ of characteristic $0$, such that
$[\D^{i},\D^j]\subset\D^{i+j-1}$, where $[\cdot,\cdot]$ is the commutator
bracket and where $\D^{i}=\{ 0\}$ for $i<0$, by convention.

It is obvious that $\A =\D^0$ is a commutative subalgebra of $\D$
(we will call it the \textit{basic algebra} of $\D$) and $\D^1$ a
Lie subalgebra of $\D$. We shall refer to elements $k$ of $K$,
naturally embedded in $\A$ or $\D$ by $k\in K\rightarrow
k1\in\A\subset\D$, as \textit{constants}, to elements $f$ of $\A$
as \textit{functions} and to elements $D$ of $\D$ as
\textit{differential operators}. One easily sees that every
element      $D\in\D^1$,      i.e.       every \textit{first-order
differential operator}, induces a derivation $\hat{D}\in Der(\A)$
of $\A$ by $\hat{D}(f)=[D,f]$.

By a \textit{classical  Poisson  algebra}  we  understand  a
commutative associative algebra with an $\N$-gradation
$S=\oplus_{i=0}^\infty  S_i$, $S_iS_j\subset S_{i+j}$, with unit 1
over a field $K$  of  characteristic 0, equipped with a  Poisson
bracket  $\{\cdot,\cdot\}$  such  that  $\{ S_i,S_j\}\subset
S_{i+j-1}$. Of course, we can think  of  $S$  as  of  a
$\Z$-graded algebra putting  $S_i=\{  0\}$  for   $i<0$, and as a
filtered algebra putting $S^i=\oplus_{k\le i}S_k$.  Like in the
case  of   the quantum Poisson algebra, $\A=S_0$ is an associative
and  Lie-commutative subalgebra of $S$  (the \textit{basic
algebra})  and  $S^1$  is  a  Lie subalgebra of
$(S,\{\cdot,\cdot\})$ acting on $\A$ by derivations.

An operator $\phi\in Hom_K(V_1,V_2)$ between $\N$-filtered vector spaces
\textit{respects the filtration}, if $\phi(V_1^{i})\subset V_2^{i}$ and
\textit{is lowering}, if $\phi(V_1^{i})\subset V_2^{i-1}$.

Quantum Poisson algebras induce canonically classical  Poisson  algebras
as follows. For a quantum Poisson algebra $\D$ consider
the graded vector space $S(\D)=\oplus_{i\in \Z}S_i(\D)$,
$S_i(\D)=\D^{i}/\D^{i-1}$. We have the obvious canonical surjective map
$\sigma:\D\rightarrow S$, the \textit{principal-symbol map}. Note that
$\sigma(\A)=\A=S_0$.  By $\sigma(D)_j$  we  denote  the
projection of $\sigma(D)$ to  $S_j$.

Since for each non-zero differential operator $D\in\D$, there is a
single $i=deg(D)\in\Z$ such that $D\in\D^{i}\setminus\D^{i-1}$,
$\sigma(D)_j=0$  if  $j\ne deg(D)$ and $\sigma(D)_{deg(D)}=\sigma(D)$.
We  set,   for   $\dot{D_1}\in   S_i$, with
$\dot{D_1}=\sigma(D_1)$,
and $\dot{D_2}\in  S_j$,  with  $\dot{D_2}=\sigma(D_2)$,
\[
\dot{D_1}\dot{D_2}=\sigma(D_1\cdot D_2)_{i+j},\quad
\{\dot{D_1},\dot{D_2}\}=\sigma([D_1,D_2])_{i+j-1}.
\]
It is easy to  see that these definitions do not depend on the choice of
the representatives $D_1$ and $D_2$ and that we get a classical  Poisson
algebra with the same basic algebra $\A$. This classical Poisson algebra
we call the \textit{classical limit  of  the  quantum  Poisson  algebra}
$\D$. We can formulate this as follows.
\begin{theo} For every quantum Poisson algebra $\D$ there is a unique
classical Poisson algebra structure on the graded vector space $S(\D)$
such that
\begin{equation}\label{Poi1}
\sigma(D_1)\sigma(D_2)=\sigma(D_1\cdot D_2)_{deg(D_1)+deg(D_2)}
\end{equation}
and
\begin{equation}\label{Poi2}
\{\sigma(D_1),\sigma(D_2)\}=\sigma([D_1,D_2])_{deg(D_1)+deg(D_2)-1}
\end{equation}
for each $D_1,D_2\in\D$. In particular,
\[
\{\sigma(D_1),\sigma(D_2)\}=\begin{cases} \sigma([D_1,D_2])&\\
\mbox{or}&\\
0.&\end{cases}\]\label{lem1}
\end{theo}

\begin{cor}
For $D_1,D_2,\ldots,D_n\in\D$, if
\[[D_1,[D_2,\ldots,[D_{n-1},D_n]]]=0,\] then
\[\{\s(D_1),\{\s(D_2),\ldots,\{\s(D_{n-1}),\s(D_n)\}\}\}=0.\]
\end{cor}
Note that every  linear  map  $\Phi:\D_1\raa\D_2$  between  two
quantum Poisson algebras, which respects the  filtration,  induces
canonically  a linear map $\tilde\Phi:S(\D_1)\raa S(\D_2)$, which
respects the  gradation, by $\tilde\Phi(\sigma(D))=\s(\Phi(D))$.
In view of Theorem \ref{lem1}, it is easy to see that if such
$\Phi$ is  a  homomorphism  of  associative (resp.  Lie)
structure, then  $\tilde\Phi$  is  a   homomorphism   of
associative (resp. Lie) structure.

A classical Poisson algebra $S$ is said  to  be
\textit{non-singular}, if   $\{   S^1,\A\}=\A$.  The   Poisson
algebra    $S$    is    called \textit{symplectic}, if constants
are the only central elements of $(S,\{\cdot,\cdot\})$, and
\textit{distinguishing}, if for any $P\in S$ one has: \[\forall
f\in\A, \exists n\in\N:
\underbrace{\{P,\{P,\ldots,\{P}_{n},f\}\}\}=0\Rightarrow P\in\A.\]
A  quantum  Poisson  algebra  is  called  \textit{non-singular}
(resp. \textit{symplectic} or \textit{distinguishing}), if its
classical  limit is a non-singular (resp. symplectic or
distinguishing) classical  Poisson algebra.
\begin{prop} For any quantum Poisson algebra $\D$:
\begin{description}
\item{(a)} $\D$ is non-singular if and only if $[\D^1,\A]=\A$;
\item{(b)} if $\D$ is symplectic, then the constants are the  only
central elements in $\D$;
\item{(c)} if $\D$ is distinguishing, then for any $D\in\D$ one
has: \[\forall
f\in\A, \exists n\in\N:
\underbrace{[D,[D,\ldots,[D}_{n},f]]]=0\Rightarrow
D\in\A.\]
\end{description}
\label{prop0}
\end{prop}
\textit{Proof.} It is obvious that
\[ [\D^1,\A]=\sigma([\D^1,\A])=\{ S^1(\D),\A\},
\]
which proves (a). To prove the part (b), it suffices to observe
that the center of the Lie algebra $S(\D)$  contains the image of
the center of $\D$ by the map $\sigma$. Finally, in view of
Corollary 1,
\[\underbrace{[D,[D,\ldots,[D}_{n},f]]]=0\]
implies
\[\underbrace{\{\sigma(D),\{\sigma(D),\ldots,\{\sigma(D)}_{n},f\}\}\}=0
\]
and part (c) follows.\rule{1.5mm}{2.5mm}

\medskip\noindent
\textbf{Example 1} A standard example of a quantum Poisson algebra
is the  algebra  $\D(M)$ of differential operators
$D:C^{\infty}(M)\rightarrow C^{\infty}(M)$ associated with a
manifold $M$. Its classical limit $S(M)$ is the Poisson algebra
$Pol(T^*M)$ of polynomials on the cotangent bundle $T^*M$ (i.e. of
the smooth functions on $T^*M$ that are polynomial along the
fibers) with the standard symplectic Poisson bracket on $T^*M$.
One can view also $S(M)$ as  the  algebra  of symmetric
contravariant tensors  on  $M$  with  the  symmetric  Schouten
bracket. We have a canonical splitting \[\D(M)=\A\oplus\D_c(M),\]
where $\A=C^{\infty}(M)$ and where $\D_c(M)$ is the algebra of
differential operators vanishing on constants ($D\in\D_c(M)$ if
and only if $D(1)=0$). If $\D^{i}_c(M)=\D^{i}(M)\cap\D_c(M)$
($i\ge 0$), we also have $\D^{i}(M)=\A\oplus\D^{i}_c(M)$. It is
clear that $\D^0_c(M)=0$ and that $\D^1_c(M)$ is the Lie algebra
$Der(\A)$ of derivations of $\A$, i.e. the Lie algebra $Vect(M)$
of vector fields on $M$. Note  that the Lie algebras $\D^1(M)$ and
$S^1(M)$  are  both  isomorphic  to $Vect(M)\oplus C^\infty(M)$
with the bracket $[X+f,Y+g]=[X,Y]+(X(g)-Y(f))$.

The quantum Poisson algebra $\D(M)$ is easily seen to be
non-singular and symplectic. We will show in the next section that it is
distinguishing.

\medskip\noindent
\textbf{Example 2} The above example can be extended to the case
of the quantum Poisson algebra of differential operators on a
given associative commutative algebra $\A$ with unit  1.  The
corresponding  differential calculus has been developed and
extensively studied by  A.~M.~Vinogradov \cite{Vi}.

\medskip
To investigate the algebra $\D(M)$ of differential operators we
need some preparations. Let us look at local representations  of
differential operators and the formal calculus (see e.g.
\cite{DWLc}, \cite{NP2}).

Consider an open subset $U$ of $\R^n$, two real finite-dimensional
vector spaces $E$ and $F$, and some local operator
\[O\in {\cal L}(C^\infty(U,E),C^\infty(U,F))_{loc}.\]
The operator is fully defined by its values on the products $fe$,
$f\in C^{\infty}(U), e\in E$. A well known theorem of J. Peetre
(see \cite{JP}) states that it has the form
\[O(fe)=\sum_{\alpha}O_{\alpha}(\partial^{\alpha}(fe))
=\sum_{\alpha}O_{\alpha}(e)\partial^{\alpha}f,\] where
$\partial_x^{\alpha}=\partial_{x^1}^{\alpha^1}\ldots
\partial_{x^m}^{\alpha^m}$
and $O_{\alpha}\in C^{\infty}(U,{\cal L}(E,F))$. Moreover, the
coefficients $O_{\alpha}$ are well determined by $O$ and the
series is locally finite (it is finite, if $U$ is relatively
compact).

We shall symbolize the partial derivative $\partial^{\alpha}f$ by
the monomial
$\xi^{\alpha}=\xi_1^{\alpha^1}\ldots\linebreak\xi_m^{\alpha^m}$ in
the components $\xi_1,\ldots,\xi_m$ of some linear form $\xi\in
(\R^n)^*$, or---at least mentally---even by $\xi^{\alpha}f,$ if
this is necessary to avoid confusion. The operator $O$ is thus
represented by the polynomial
\[{\cal O}(\xi;e)=\sum_{\alpha}O_{\alpha}(e)\xi^{\alpha}.\]
When identifying the space $Pol((\R^n)^*)$ of polynomials on
$(\R^n)^*$ with the space $\vee \R^n$ of symmetric contravariant
tensors of $\R^n$, one has ${\cal O}\in C^{\infty}(U,\vee
\R^n\otimes {\cal L}(E,F))$. Let us emphasize that the form $\xi$
symbolizes the derivatives in $O$ that act on the argument $fe\in
C^{\infty}(U,E)$, while $e\in E$ represents this argument. In the
sequel, we shall no longer use different notations for the
operator $O$ and its representative polynomial ${\cal O}$; in
order to simplify notations, it is helpful to use even the same
typographical sign, when referring to the argument $fe$ and its
representation $e$.

Let us for instance look for the local representation of the Lie
derivative of a differential operator (it is well-known that
$L_XD=[X,D]$ ($X\in Vect(M)$, $D\in\D^{i}(M)$ or
$D\in\D^{i}_c(M)$) defines a module structure over $Vect(M)$ on
$\D^{i}(M)$ resp. $\D^{i}_c(M)$). If $D\in {\D(M)}$, its
restriction $D\!\!\mid_U$ (or simply $D,$ if no confusion is
possible) to a domain $U$ of local coordinates of $M$, is a local
operator from $C^\infty (U)$ into $C^\infty(U)$ that is
represented by $D(f)\simeq D(\xi;1)=D(\xi)$, where $f\in C^\infty
(U)$ and where $\xi$ represents the derivatives acting on $f$. The
Lie derivative of $D(f)$ with respect to a vector field $X\in
C^\infty (U,\R^n)$, is then represented by $L_X(D(f))\simeq \la
X,\eta+\xi\ra D(\xi)$. Here $\eta\in (\R^n)^*$ is associated to
$D$ and $\la X,\eta+\xi\ra$ denotes the evaluation of $X\in \R^n$
on $\eta +\xi$. When associating $\zeta$ to $X$, one gets
$D(L_Xf)\simeq \la X,\xi\ra D(\xi+\zeta)$ and
\begin{equation}
(L_X D)(f) \simeq \la X,\eta\ra D(\xi) - \la X,\xi\ra \tau_\zeta
D(\xi),\end{equation} where $\tau_\zeta
D(\xi)=D(\xi+\zeta)-D(\xi)$.

\section{Algebraic characterization of a manifold}

\begin{theo}  The  quantum  Poisson  algebra  $\D(M)$  of   differential
operators on $C^\infty(M)$ is distinguishing (i.e. the classical
Poisson algebra $S(M)$ is distinguishing).\end{theo}

\textit{Proof.} Since for $P,Q\in S=Pol(T^*M)$, $\{P,Q\}=H_P.Q$,
where $H_P$ is the Hamiltonian vector field of $P$, we have to
prove that if $P\in S\setminus\A$ ($\A=\ci(M)$), there is a
function $f\in\A$ such that for every integer $n\in \N$,
$(H_P)^n.f\neq 0$.

If $(U,(x^1,\ldots,x^n))$ is a chart of $M$, then $H_P$ has in the
associated Darboux chart
$(T^*U,(x^1,\ldots,x^n,\xi_1,\ldots,\xi_n))$ the classical
expression $H_P=\overline{\partial}_i P\partial_i-\partial_i
P\overline{\partial}_i$, where
$\overline{\partial}_i=\partial/\partial \xi_i$ and
$\partial_i=\partial/\partial x^{i}$. It follows from the
hypothesis $P\in S\setminus\A$ that $\overline{\partial}_i
P\partial_i\neq 0$ for at least one chart $U$ of $M$. In order to
simplify notations, we shall write in the associated Darboux chart
$T^*U$,
\[H_P=F^{i}\partial_i+G_i\overline{\partial}_i,\] with $F^{i}\partial_i\neq
0$.

First notice that, for an arbitrary neighborhood $]a,b[$ of an
arbitrary point $x_0\in \R$, it is possible to construct a
sequence $x_1,x_2,\ldots\in ]a,b[$ with limit $x_0$ and a function
$h\in C^{\infty}(\R)$ such that, if $d_x^kh$ denotes the $k$-th
derivative of $h$,
\[(d_x^kh)(x_n)\left\{\begin{array}{l}=0,\mbox{  for all }k\in
\{0,\ldots,n-1\}\\\mbox{and}\\\neq 0,\mbox{  for
}k=n\end{array}\right..\] Indeed, set $d=\frac{b-x_0}{2}$,
$x_n=x_0+\frac{d}{n}$ $(n\in \N^*)$, $\delta_n=x_n-x_{n+1}$, and
$V_n=]x_n-\frac{\delta_n}{2},x_n+\frac{\delta_n}{2}[$. It is clear
that the intersections $V_n\cap V_{n+1}$ are empty. Take now
smooth functions $\alpha_n$ with value $1$ around $x_n$ and
compact support in $V_n$ and define smooth functions $h_n$ by
$h_n(x)=(x-x_n)^n\alpha_n(x)$. One easily sees that
$(d_x^kh_n)(x_n)$ vanishes for all $k\in\{0,\ldots,n-1\}$ and does
not for $k=n$. Finally, the function $h$ defined by
$h(x)=\sum_{n=1}^{\infty}h_n(x)$ has all the desired properties.

When returning to the initial problem, remark that at least one
$F^{i}$ does not vanish, say $F^1$. If its value at some point
$(x_0,\xi_0)\in T^*U$ is non-zero, the function
$F^1(\cdot,\xi_0)\in C^{\infty}(U)$ is non-zero on some
neighborhood $V$ of $x_0$.

In the sequel, the coordinates $(x^1,\ldots,x^n)$ of a point $x\in
U$ will be denoted by $(x^1,x'')\in \R\times \R^{n-1}$. Consider
now $V$ as an open subset of $\R^n$, introduce the section
$V^1=\{x^1:(x^1,x_0'')\in V\}$ of $V$ at the level $x_0''$, and
construct the previously described sequence $x^1_n$ and function
$h$ in this neighborhood $V^1$ of $x_0^1$. The sequence defines a
sequence $x_n$ in $V$ with limit $x_0$ and the function defines a
function still denoted by $h$ in $C^{\infty}(V)$.

When multiplying this $h$ by a smooth $\alpha$, which has value
$1$ in a neighborhood of the $x_n$'s and is compactly supported in
$V$, we get the function $f\in C^{\infty}(M)$ that we have to
construct. Indeed, for every $n$,
\[\left((H_P)^n.f\right)(x_n,\xi_0)=\left((F^{i}\partial_i+G_i
\overline{\partial}_i)^nh\right)(x_n,\xi_0).\]
The function on the r.h.s. is a sum of terms in the
$\p_ih,\p_{i_1}\p_{i_2}h,\ldots,\p_{i_1}\ldots\p_{i_n}h$ and the
maximal order terms are $F^{i_1}\ldots
F^{i_n}\p_{i_1}\ldots\p_{i_n}h$. All the terms of order less than
$n$ vanish, since the derivatives with respect to $x^{i}$ ($i\neq
1$) vanish and for $k<n$, $(d_{x^1}^kh)(x^1_n)=0$. The terms of
maximal order $n$ also vanish, except $(F^1)^nd_{x^1}^nh$ that is
non-zero at $(x_n,\xi_0)$.\rule{1.5mm}{2.5mm}

\medskip
For any Lie algebra $(\Li,[\cdot,\cdot])$, by $Nil(\Li)$ we denote
the set of those $D\in\Li$ for which $ad_D$ is locally nilpotent:
\[ Nil(\Li)=\{D\in\Li:\forall D'\in\Li,\exists n\in
\N:\underbrace{[D,[D,\ldots,[D}_{n},D']]]=0\}.
\]
\begin{prop}
If a quantum or classical Poisson algebra $\Li$ with the basic
algebra $\A$ is distinguishing, then
\begin{description}
\item{(a)} $Nil(\Li)=\A$,
\item{(b)}
\[\{ P\in S:\{ P,\A\}\subset S_{i}\}=S_{i+1}\oplus\A, \quad (i\ge -1)
\]
in case $\Li=S$ is classical. In particular,
\[\{ P\in S:\{ P,\A\}\subset S^{i}\}=S^{i+1},\quad(i\ge -1).
\]
\item{(c)}
\[\{ D\in\D:[D,\A]\subset\D^{i}\}=\D^{i+1},\quad(i\ge -1)
\]
in case $\Li=\D$ is quantum.
\end{description}
\label{prop1}\end{prop}
\textit{Proof.} (a) is obvious for classical and, in view  of Proposition
\ref{prop0}, also for quantum Poisson algebras. \\
(b)  Since  $\{ P,\A\}\subset S\ominus S_i$,  for  any  $P\in
S\ominus
S_{i+1}$, the inclusion $\{  P,\A\}\subset  S_i$ for such $P$ implies
$\{ P,\A\}=0$,   so $P\in\A$.\\
(c) If $D\in\D\setminus\D^{i+1}$ and $[D,\A]\subset\D^{i}$ ($i\ge
-1$), then $\{\s(D),\A\}=0$ and
$\s(D)\in\A$, which is contradictory.\rule{1.5mm}{2.5mm}

\medskip
Now we will start the studies on properties of  isomorphisms  of  quantum
and classical Poisson algebras. We will concentrate on the quantum level,
since on the classical level all the considerations  are  analogous  and
even simpler.

\begin{cor}
Every  isomorphism  $\Phi:\D_1\raa\D_2$   of    the
Lie algebras $(\D_i,[\cdot,\cdot])$  for   distinguishing  quantum
Poisson  algebras $\D_i$, $i=1,2$,
respects    the    filtration    and     induces     an isomorphism
$\tilde{\Phi}:S(\D_1)\raa S(\D_2)$,
$\tilde{\Phi}(\sigma(D))=\sigma(\Phi(D))$,    of
the corresponding classical limit Lie algebras.\label{cor2}
\end{cor}
\textit{Proof.} It is obvious that $\Phi(Nil(\D_1))=Nil(\D_2)$,
i.e. $\Phi(\A_1)=\A_2$.
Inductively, if $\Phi(\D_1^{i})\subset\D_2^{i}$, then, for any
$D\in\D_1^{i+1}$, \[[\Phi(D),\A_2]=\Phi([D,\A_1])\subset\D_2^{i},\]and
$\Phi(D)\in\D_2^{i+1}$, by Proposition \ref{prop1}. Now, since $\Phi$  and
$\Phi^{-1}$  respect  the  filtration,   $\tilde{\Phi}$   is   a   linear
isomorphism of $S_1$ onto $S_2$ which, as easily seen,  is  a  Lie
algebra isomorphism.\rule{1.5mm}{2.5mm}\\

Denote by $C(\D)$ the centralizer of $ad_{\A}$ in $Hom_K(\D,\D)$:
\[\Psi\in C(\D)\Leftrightarrow
[\Psi,ad_{\A}]=0.\] Note that multiplications $m_f:\D\ni
D\rightarrow f\cdot D\in\D$ and
$m'_f:\D\ni D\rightarrow D\cdot f\in\D$
by elements $f\in\A$, belong to $C(\D)$.

\begin{theo}
Assume that $\D$ is a non-singular and distinguishing quantum
Poisson algebra. Then any $\Psi\in C(\D)$ respects the filtration
and there is an $f\in\A$ and a lowering $\Psi_1\in C(\D)$, such
that
\[\Psi=m_f+\Psi_1.\]\label{theo2}
\end{theo}
\textit{Proof.} (i) Information $[\Psi,ad_{\A}]=0$ means that
\begin{equation}[\Psi(D),f]=\Psi([D,f]),\label{centr1}\end{equation}
for all $D\in\D$ and
all $f\in\A$. For $D\in\A$ we get $[\Psi(D),f]=0$, so
$\Psi(D)\in\A$. Inductively, if $\Psi(\D^{i})\subset\D^{i}$, then
(\ref{centr1}) implies $[\Psi(\D^{i+1}),f]\subset\D^{i}$ and
$\Psi(\D^{i+1})\subset\D^{i+1}$.

(ii) Let now $D\in\D^1$: $\Psi(D)\in\D^1$. Since for any $f\in\A$
\[2\Psi(f[D,f])=\Psi([D,f^2])=[\Psi(D),f^2]=2f[\Psi(D),f]=2f\Psi([D,f]),\]
we                                                                  have
\begin{equation}\Psi(f\hat{D}(f))=f\Psi(\hat{D}(f)),\label{centr2}
\end{equation}
for any $f\in\A,D\in\D^1$. Substituting $D:=gD$
$(g\in\A,D\in\D^1)$ and $f:=f+h$ $(f,h\in\A)$ in (\ref{centr2}),
we get
\begin{equation}\Psi(fg\hat{D}(h))+\Psi(gh\hat{D}(f))=f\Psi(g\hat{D}(h))+
h\Psi(g\hat{D}(f)).\label{centr3}\end{equation}
For $g=\hat{D}(h)$, equation (\ref{centr3}) reads
\[\Psi(f(\hat{D}(h))^2)+\Psi(h\hat{D}(f)\hat{D}(h))=
f\Psi((\hat{D}(h))^2)+h\Psi(\hat{D}(f)\hat{D}(h)),\]
where the last terms of the l.h.s. and the r.h.s. cancel in view
of (\ref{centr2}) applied to $D:=\hat{D}(f)D$. Hence for each
$f,h\in\A$ and $D\in\D^1$,
\[\Psi(f(\hat{D}(h))^2)=f\Psi((\hat{D}(h))^2).\] The last
equation shows that the radical $rad(J)$ of the ideal
$J=\{g\in\A:\Psi(fg)=f\Psi(g),\forall f\in\A\}$ of the associative
commutative algebra $\A$, contains $[\D^1,\A]$. Since $\D$ is
non-singular, this implies that $J=\A$, so that
\[\Psi(f)=\Psi(1)f,\] for all $f\in\A$. It is obvious that
\[\Psi_1=\Psi-m_{\Psi(1)}\] belongs to $C(\D)$ and respects the
filtration, and one easily sees that it is lowering. Indeed, since
$\Psi_1(\A)=0$, assume inductively that
$\Psi_1(\D^{i})\subset\D^{i-1}$. Then,
$[\Psi_1(\D^{i+1}),\A]=\Psi_1([\D^{i+1},\A])\subset\D^{i-1}$ and
$\Psi_1(\D^{i+1})\subset\D^{i}$.\rule{1.5mm}{2.5mm}

\begin{theo}
Let $\D_i$ be distinguishing, non-singular and symplectic, $i=1,2$. Then
every isomorphism $\Phi:\D_1\raa\D_2$ of  the  Lie  algebras
$(\D_i,[\cdot,\cdot])$, $i=1,2$,  respects  the   filtration   and   its
restriction $\Phi\mid_{\A_1}$ to $\A_1$ has the form
\[\Phi\mid_{\A_1}=\kappa A,\]
where $\kappa\in K, \kappa\neq 0$ and $A:\A_1\raa \A_2$ is an
isomorphism of the associative commutative algebras. The  same  is
true  for any isomorphism  $\Phi:\D_1^1\raa\D_2^1$  of   the
corresponding Lie algebras of first-order differential operators.
\label{theo3}
\end{theo}
\textit{Proof.} By Corollary \ref{cor2}, $\Phi$ respects the
filtration, so $\Phi(\A_1)=\A_2$. Let
\[\Phi_{*}:Hom_K(\D_1,\D_1)\raa Hom_K(\D_2,\D_2)
\]
be the induced isomorphism of the Lie algebras of linear homomorphisms,
defined for $\Psi\in Hom_K(\D_1,\D_1)$ by
:\[\Phi_{*}(\Psi)=\Phi\circ\Psi\circ\Phi^{-1}.\]
Since $\Phi(\A_1)=\A_2$, $\Phi_{*}(C(\D_1))=C(\D_2)$; in particular
$\Phi_{*}(m_g)\in C(\D_2)$  for $g\in\A_1$. By Theorem \ref{theo2},
\[\Phi_{*}(m_g)(f')=\Phi_{*}(m_g)(1)\cdot f',\] i.e.
\begin{equation}\Phi(g\cdot\Phi^{-1}(f'))=\Phi(g\cdot\Phi^{-1}(1))\cdot
f',\label{aut1}\end{equation} for all $f'\in\A_2$. Observe that
$\Phi^{-1}(1)$ is central in $\D_1$ and is thus a non-vanishing
constant $\kappa^{-1}$. Substituting $\Phi(f)$ ($f\in\A_1$) to $f'$ in
(\ref{aut1}), one obtains
\[\Phi(f\cdot g)=\kappa^{-1}\Phi(f)\cdot\Phi(g).\] For $A$ defined by
\[A(f)=\kappa^{-1}\Phi(f),\] this reads \[A(f\cdot g)=A(f)\cdot A(g),\] which
completes the proof of Theorem \ref{theo3}.\rule{1.5mm}{2.5mm}\\

We can prove in the same way---mutatis mutandis---that Theorem
\ref{theo3} is still valid for $\D_i$, $i=1,2$,  replaced  by  classical
Poisson algebras $S_i$, $i=1,2$.

\begin{theo}
Let $S_i$ be a distinguishing,  non-singular  and  symplectic  classical
Poisson algebra, $i=1,2$. Then
every isomorphism $\Phi:S_1\raa S_2$ of  the  Lie  algebras
$(S_i,\{\cdot,\cdot\})$, $i=1,2$,  respects  the   filtration   and   its
restriction $\Phi\mid_{\A_1}$ to $\A_1$ has the form
\[\Phi\mid_{\A_1}=\kappa A,\]
where $\kappa\in K, \kappa\neq 0$ and $A:\A_1\raa \A_2$ is an
isomorphism of the associative commutative algebras.
\label{theo40}
\end{theo}

\begin{cor}  If  two   distinguishing,   non-singular   and   symplectic
quantum  (resp.  classical)  Poisson  algebras  are  isomorphic  as  Lie
algebras, then  their  basic  algebras  are  isomorphic  as
associative  algebras.  The  same  remains  true  for  Lie subalgebras  of
first-order operators of such Poisson algebras: if they  are  isomorphic,
then their basic algebras are isomorphic associative algebras.
\end{cor}

Let us now return to the quantum Poisson algebra $\D=\D(M)$ of
differential operators of a smooth, Hausdorff, second countable,
connected manifold $M$.
It is well known that every associative algebra isomorphism
$A:\A_1=C^{\infty}(M_1)\rightarrow\A_2=C^{\infty}(M_2)$ is of the form
\[A:\A_1\ni f\rightarrow f\circ\phi^{-1}\in\A_2,\] where
$\phi:M_1\raa M_2$ is a diffeomorphism.
Thus, we can draw a conclusion of the same type than a classical result of
Pursell and Shanks \cite{PS},\cite{JG}:
\begin{theo}The Lie algebras $\D(M_1)$ and $\D(M_2)$ $($resp. $\D^1(M_1)$
and $\D^1(M_2)$, or $S(M_1)$ and $S(M_2)$$)$ of all differential operators
$($resp. all differential  operators  of  order  $1$,  or  all  symmetric
contravariant tensors$)$ on two smooth manifolds
$M_1$ and $M_2$ are isomorphic if and only if the manifolds $M_1$
and $M_2$ are diffeomorphic.
\end{theo}
Studying the isomorphisms mentioned in the above theorem reduces then  to
studying automorphisms  of  the  Lie  algebras  $\D(M)$,  $\D^1(M)$,  and
$S(M)$.

\section{Particular automorphisms}
In the sequel, $\D$ denotes the quantum algebra $\D(M)$ and $S$ is
its classical limit $S(M)$.

\textbf{1.}  Every  automorphism  $A$ of   the   associative
algebra $\A=C^\infty(M)$ (which is implemented by a diffeomorphism
$\phi$ of $M$) induces an automorphism $A_{*}$ of the Lie algebra
$\D$:\[A_{*}(D)=A\circ D\circ A^{-1}\]($D\in\D$). It clearly
restricts to an automorphism of $\D^1$. The automorphism $A$
induces also an automorphism $A_*$ of $S$. It is just induced by
the  phase  lift of  the  diffeomorphism $\phi$ to the cotangent
bundle  $T^*M$  if we  interpret  elements  of $S$ as polynomial
functions on $T^*M$. If  we  interpret  $S$  as symmetric
contravariant tensors on $M$, then $A_*$ is just the action of
$\phi$  on such tensors.

Let    now    $\Phi\in    Aut(\D,[\cdot,\cdot])$    (resp.
$\Phi\in Aut(\D^1,[\cdot,\cdot])$ or  $\Phi\in
Aut(S,\{\cdot,\cdot\})$).   By Theorem \ref{theo3}, there are
$A\in Aut(\A,\cdot)$ and $\kappa\in K, \kappa\neq 0$, such that
\[\Phi\!\!\mid_{\A}=\kappa A_*\!\!\mid_{\A}.\] Then,
\[\Phi_1=(A_*)^{-1}\circ\Phi\]is an automorphism of $\D$ (resp.
$\D^1$ or $S$), which is $\kappa\!\cdot\! id$ ($id$ is the identity map) on
$\A$. It is thus sufficient to describe the automorphisms that are
$\kappa\!\cdot\! id$ on functions.

\textbf{2.} Let $\omega\in\Omega^1(M)\cap\ker d$ be a closed
$1$-form on $M$ and $D\in\D^{i}$. If $U$ is an open subset of $M$
and $\omega\!\!\mid_U=d(f_U)$ ($f_U\in C^{\infty}(U)$), the
operators $[D\!\!\mid_U,f_U]\in\D^{i-1}_U$ ($\D^k_U$ is defined as
$\D^k$ but for $M=U$) are of course the restrictions of an unique
well-defined operator
$\overline{\omega}(D)\in\D^{i-1}$:\[\overline{\omega}
(D)\!\!\mid_U=[D\!\!\mid_U,f_U],\] since the above commutator does
not depend on the choice of $f_U$ with $\omega\!\!\mid_U=d(f_U)$
(constants are  central  with  respect  to  the bracket). It is
clear that $\overline{\omega}\in {\cal L}(\D,\D)\cap {\cal
L}(\D^{i},\D^{i-1})$, that $\overline{\omega}(X)=\omega(X)$ for
all $X\in Vect(M)$, and that $\omega\mapsto\overline{\omega}$ is
linear. Moreover, $\overline{\omega}$ is a $1$-cocycle of the
adjoint Chevalley-Eilenberg cohomology of $\D$, i.e. a derivation
of $\D$. Since $\overline{\omega}$ is lowering, it is locally
nilpotent, so that
\[e^{\overline{\omega}}=id+\overline{\omega}+\frac{1}{2!}
\overline{\omega}^2+\ldots\] is well defined and it is  an
automorphism  of  $\D$  (that  is  identity  on  functions).  In
particular, for $\omega=df$, the automorphism
$e^{\overline{\omega}}$ is just the inner automorphism $\D\ni
D\mapsto e^f\cdot D\cdot e^{-f}\in\D$.

On the classical level we have an analogous derivation of the  classical
Poisson algebra $S$:
\[\overline{\omega}
(P)\!\!\mid_U=\{ P\!\!\mid_U,f_U\},\] and the analogous
automorphism $e^{\overline{\omega}}$. These automorphisms have a
geometric description, if we interpret $S$ as the Lie algebra of
polynomial functions  on  the cotangent bundle $T^*M$ with the
canonical Poisson bracket. Every  closed  1-form  $\omega$  on $M$
induces  a  vertical   locally hamiltonian  vector  field
$\omega^v$  on  $T^*M$ which  connects  the 0-section of $T^*M$
with another lagrangian  submanifold which  is  the image of the
section $\omega$. If, locally, $\omega=df$, then $\omega^v$  is,
locally,  the  hamiltonian  vector  field  of  the pull-back of
$f$ to $T^*M$. In the pure vector bundle language, $\omega^v$ is
just the vertical lift of the section $\omega$ of  $T^*M$.  Since
this vector field is vertical and constant on  fibers,  it  is
complete  and determines a one-parameter group $Exp(\omega^v)$  of
symplectomorphisms of $T^*M$. The automorphism
$e^{\overline{\omega}}$ is just  the  action of   $Exp(\omega^v)$
on   polynomial   functions   on   $T^*M$.    The
symplectomorphism $Exp(\omega^v)$ translates every covector
$\eta_p$  to
$\eta_p+\omega(p)$.\\

\textbf{3.} The following remark concerns the divergence operator
on an arbitrary manifold $M$. For further details the reader is
referred to \cite{PBAL}.

Denote by $I\!\!F_{\lambda}(TM)$ ($\lambda\in\R$) the vector
bundle (of rank $1$) of $\lambda$-densities and by
$\F_{\lambda}(M)$ the $Vect(M)$-module of $\lambda$-density fields
(or simply $\lambda$-densities) on $M$ (i.e. the space of smooth
sections of $I\!\!F_{\lambda}(TM)$, endowed with the natural Lie
derivative $L_X,$ $X\in Vect(M)$). The result stating that these
modules $\F_{\lambda}(M)$ are not isomorphic, implies the existence
of a non-trivial $1$-cocycle of the Lie algebra $Vect(M)$
canonically represented on $\ci(M)$. It appears, if
$\F_{\lambda}(M)$ is viewed as a deformation of $\F_0(M)=\ci(M)$.

Let us be somewhat more precise. In the proof of triviality of the
bundles $I\!\!F_{\lambda}(TM)$, one constructs a section that is
everywhere non-zero (and even, which has at each point only
strictly positive values). Let $\rho_0\in\F_1(M)$ be such a
section. Then $\rho_0^{\lambda}\in\F_{\lambda}(M)$ also vanishes
nowhere and $\tau_0^{\lambda}:f\in\ci(M)\longrightarrow
f\rho_0^{\lambda}\in\F_{\lambda}(M)$ is a bijection. One has the
subsequent results:
\begin{itemize}
\item There is a $1$-cocycle $\gamma:Vect(M)\longrightarrow\ci(M)$,
which depends on $\rho_0$ but not on $\lambda$, such that, for any
$X\in Vect(M)$,
\[(\tau_0^{\lambda})^{-1}\circ L_X\circ\tau_0^{\lambda}:f\in\ci(M)\longrightarrow
X(f)+\lambda\gamma(X)f\in\ci(M).\] \item The cocycle $\gamma$ is a
differential operator with symbol $\sigma(\gamma)(\zeta;X)=\la
X,\zeta\ra,$ where $\la X,\zeta\ra$ denotes the evaluation of
$\zeta\in T_x^*M$ upon $X\in T_xM$.\item The cohomology class of
$\gamma$ is independent of $\rho_0$.
\end{itemize}
This class $div_M$ is the class of the divergence. Each cocycle
cohomologous to $\gamma$ will be called a divergence. Finally, the
following propositions hold: \begin{itemize}\item The first
cohomology space of $Vect(M)$ represented upon $\ci(M)$ is given
by \begin{equation} H^1(Vect(M),\ci(M))=\R\;div_M\oplus
H^1_{DR}(M),\label{cohovect}\end{equation} where $H^1_{DR}(M)$
denotes the first space of the de Rham cohomology of $M$. \item
For any divergence $\gamma$ on $M$, there is an atlas of $M$, such
that in every chart,
\begin{equation}\gamma(X)=\sum_i\p_{x^{i}}X^{i}, \forall X\in Vect(M),\label{classdiv}\end{equation} with
self-explaining notations.
\end{itemize}

The preceding results have a simple explanation. Remember that if
the manifold $M$ is orientable and if $\Omega$ is a fixed volume
form, the divergence of $X\in Vect(M)$ with respect to $\Omega$ is
defined as the smooth function $div_{\Omega}X$ of $M$ that
verifies $L_X\Omega=(div_{\Omega}X)\mbox{ }\Omega$. One easily
sees that \[div_{-\Omega}X=div_{\Omega}X.\] But this means that
the divergence of a vector field can even be defined on a
non-orientable manifold with respect to a pseudo-volume form.

The divergence operator associated to a $1$-density $\rho_0,$ will
be denoted by $div_{\rho_0}$ or simply $div,$ if no confusion is
possible. Let us fix a divergence on $Vect(M).$

\begin{lem} There is a unique ${\cal  C}\in  Aut(\D,[\cdot,\cdot]),$  such
that ${\cal C}(f)=-f,$ ${\cal C}(X)=X+divX,$
\[{\cal C}(D\circ f)=f\circ {\cal C}(D),\] and
\[{\cal C}(D\circ X)=-{\cal C}(X)\circ {\cal C}(D),\]
for all $f\in \A,X\in\D^1_c,$ and $D\in\D$. \label{C}\end{lem}

\textit{Proof.}  Consider an atlas of $M,$ such that the
divergence has the form (\ref{classdiv}) in any chart. Then, in
every chart $(U,(x^1,\ldots,x^n)),$ ${\cal C}$ given by
\[{\cal C}(\eta;P_k)(\xi)=(-1)^{k+1}P_k(\xi +\eta),\]
where $P_k\in\vee ^k\R^n$ is a homogeneous polynomial of degree
$k$, defines an operator ${\cal C}_U:\D_U\longrightarrow\D_U$ that
(maps $\D^{i}_U$ into $\D^{i}_U$ and) verifies the above
characteristic properties. Let's explain for instance the fourth
one, the third is analogous and the first and second are obvious.
Use the previously mentioned simplifications of notations,
identify the space $\D^{i}_U$ of differential operators to the
space $\ci(U,\vee^{\le i}\R^n)$ of polynomial representations, set
$X=gX$ and $D=hP_k$ (on the l.h.s. $X\in C^{\infty}(U,\R^n)$ and
$D\in\ci(U,\vee^{\le i}\R^n)$, on the r.h.s. $g,h\in \ci(U)$,
$X\in \R^n$, and $P_k\in\vee^k\R^n$ ($k\le i$)), and symbolize the
derivatives acting on $g,h$ and the argument $f\in\ci(U)$ of
$D\circ X$, $\C_U(D\circ X)$, and $\C_U(X)\circ\C_U(D)$, by
$\zeta,\eta$ resp. $\xi$. Since
\[(D\circ X)(f)=D(X(f))\simeq \la X,\xi\ra P_k(\xi+\zeta)=\la
X,\xi\ra \sum_{\E}\frac{1}{\E !}(\zeta\p_{\xi})^{\E}P_k(\xi)\]
($\zeta\p_{\xi}$: derivative with respect to $\xi$ in the
direction of $\zeta$), one has
\[\begin{split}\left(\C_U(D\circ X)\right)(f)&\simeq
\sum_{\E}\frac{1}{\E
!}\C(\eta+\zeta;X(\zeta\p_{\xi})^{\E}P_k)(\xi)\\
 &=(-1)^k\la X,\xi+\eta+\zeta\ra\left(\sum_{\E}\frac{1}{\E!}
((-\zeta)\p_{\xi})^{\E}P_k\right)(\xi+\eta+\zeta)\\
 &=(-1)^k\la X,\xi+\eta+\zeta\ra P_k(\xi+\eta)\\
 &\simeq -\left(\C_U(X)\circ\C_U(D)\right)(f).\end{split}\]

It is well known that any differential operator $D\in\D$ has a
global (not necessarily unique) decomposition as a finite sum of
terms of the type $f\mbox{ }{X_k}\circ\ldots\circ{X_1}$
($f\in\ci(M),X_{\E}\in Vect(M)$). If we set ${\cal
L}^1_X=X+divX$ ($X\in Vect(M)$), we have
\[\C_U(X\m_U)=\cl_X^1\m_U\] and
\[\C_U(D\m_U)=\left((-1)^{k+1}\cl^1_{X_1}\circ
\ldots\circ\cl^1_{X_k}\circ f\right)\m_U.\]
This means that the $\C_U$ are the restrictions of a unique
well-defined operator \[\C:\D\ni D=f\mbox{ }
{X_k}\circ\ldots\circ
{X_1}\longrightarrow\C(D)=(-1)^{k+1}\cl^1_{X_1}\circ
\ldots\circ\cl^1_{X_k}\circ f\in\D,\]
which inherits the characteristic properties.

The homomorphism-property, $\C[D,\Delta]=[\C D,\C\Delta]$
($D,\Delta\in\D$), is a direct consequence of the characteristic
properties and the definition of $\C$. Noticing that $\C^2(X)=X$
and---from the preceding verification---that
$\C(D\circ\Delta)=-\C(\Delta)\circ\C(D)$, one immediately sees
that $\C^2=id$, so that $\C\in Aut(\D)$.\rule{1.5mm}{2.5mm}

\begin{re} One easily convinces oneself that, if $\Omega$ is a volume
form of $M,$
$\C$ is the opposite of the conjugation $*:\D\ni D\longrightarrow
D^*\in \D$ of differential operators, defined by
\[\int_MD(f).g\mid\Omega\mid=\int_Mf.D^*(g)\mid\Omega\mid,\]
for all compactly supported $f,g\in\ci(M)$.
\end{re}

\textbf{4.} On $S$, like on every graded algebra, there is  a
canonical one-parameter  family  of   automorphisms   $U_\kappa$,
$\kappa\ne   0$,   namely $U_\kappa(P)=\kappa^{1-i}P$ for $P\in
S_i$. It is easy to see  that  $U_\kappa$  is  an automorphism of
the Lie algebra $S$. For  positive  $\kappa$  this  is  the
one-parameter group of automorphisms induced by the  canonical
derivation $Deg:S\raa S$ of  the  Poisson bracket, $Deg(P)=(i-1)P$
for $P\in S_i$, namely $U_\kappa=e^{-\log(\kappa)Deg}$. Since
$U_\kappa\!\!\mid_\A=\kappa\!\cdot\! id\!\mid_\A$,  we can   now   reduce
every automorphism  $\Phi$ of  the  Lie algebra  $S$  to   the
case when $\Phi\!\!\mid_\A=id\!\mid_\A$.

\section{Automorphisms of the Lie algebra \boldmath{$\D^1(M)$}}

When using the decomposition $\D=\A\oplus\D_c$, we denote by
$\pi_0$ and $\pi_c$ the projections onto $\A$ resp. $\D_c$.
Furthermore, if $D\in\D$, we set $D_0=\pi_0D=D(1)$ and
$D_c=\pi_cD=D-D(1)$, and if $\Phi\in {\cal L}(\D)$, we set
$\Phi_0=\pi_0\circ\Phi\in {\cal L}(\D,\A)$ and
$\Phi_c=\pi_c\circ\Phi\in {\cal L}(\D,\D_c)$. Note also that for
$f,g\in\A$, one has $[D_c,f]_0=D_c(f)$ and
$[D_c,f]_c(g)=D_c(f\cdot g)-D_c(f)\cdot g-f\cdot D_c(g)$, so that
$[D_c,f]_0=0, \forall f\in\A$ if and only if $D_c=0$ and
$[D_c,f]_c(g)=0,
\forall f,g\in\A$ if and only if $D_c\in\D^1_c$.\\

Let us now return to the problem of the determination of all
automorphisms $\Phi$ of $\D$ (resp. $\D^1$) that coincide with
$\kappa\!\cdot\! id$ on functions.\\

The projection of the homomorphism-property, written for
$D_c\in\D_c$ and $f\in\A$, leads to the equations \begin{equation}
\left(\Phi_cD_c\right)(f)=\kappa^{-1}\Phi_0[D_c,f]=D_c(f)+\kappa^{-1}\Phi_0
[D_c,f]_c\label{H1}\end{equation} and \begin{equation}
[\Phi_cD_c,f]_c=\kappa^{-1}\Phi_c[D_c,f]_c\label{H2},\end{equation}
and its projection,   if   it is   written   for
$D_c,\Delta_c\in\D_c$, gives
\begin{equation}
\Phi_0[D_c,\Delta_c]=(\Phi_cD_c)(\Phi_0\Delta_c)-(\Phi_c\Delta_c)
(\Phi_0D_c)\label{H3}\end{equation}
and
\begin{equation}\Phi_c[D_c,\Delta_c]=[\Phi_cD_c,\Phi_0\Delta_c]_c+[\Phi_0D_c,
\Phi_c\Delta_c]_c+[\Phi_cD_c,\Phi_c\Delta_c].\label{H4}\end{equation}

If one writes these equations for $D_c$ and $\Delta_c$ in the Lie
subalgebra $\D_c^1$ of $\D_c$, (\ref{H1}) means that
$\Phi_c\!\!\mid_{\D_c^1}=id$, (\ref{H2}) and (\ref{H4}) are
trivial, and (\ref{H3}) tells that
$\alpha:=\Phi_0\!\!\mid_{\D^1_c}$ is a $1$-cocycle of the Lie
algebra of vector fields canonically represented on functions by
Lie derivative. As, in view of (\ref{cohovect}), \[\alpha =\lambda
div+\omega\mbox{ }(\lambda\in \R,\omega\in\Omega^1(M)\cap\ker
d),\] one has the following

\begin{theo} The automorphisms $\Phi_1$ of $\D^1(M)$ that verify
$\Phi_1\!\!\mid_{C^{\infty}(M)}=\kappa\!\cdot\! id$ ($\kappa\in\R,
\kappa\neq 0$), are the mappings
\[\Phi_1=\kappa\pi_0+(id+\lambda div+\omega)\circ\pi_c,\]
where $\lambda\in \R$ and $\omega\in\Omega^1(M)\cap\ker d$.
\label{Phi1}\end{theo} Indeed, one easily sees that these
homomorphisms of $\D^1$ are bijective. We  can  now  summarize all
facts  and give  the complete description of automorphisms of
$\D^1$.
\begin{theo} A linear map  $\Phi:\D^1(M)\raa\D^1(M)$ is an  automorphism
of   the   Lie algebra $\D^1(M)=Vect(M)\oplus C^\infty(M)$  of
linear  first-order  differential  operators  on  $C^\infty(M)$
if and only if it can be written in the form
\begin{equation}
\Phi(X+f)=\phi_*(X)+  (\kappa f+\lambda  divX  +\omega(X)) \circ
\phi^{-1},
\end{equation}
where $\phi$ is a diffeomorphism of  $M$,  $\lambda,  \kappa$  are
constants, $\kappa\ne 0$, $\omega$ is a closed 1-form on $M$, and
$\phi_*$ is defined by
\[(\phi_*(X))(f)=(X(f\circ\phi))\circ\phi^{-1}.\] All the objects $\phi,
\lambda,\kappa,\omega$ are uniquely determined by $\Phi$.
\label{d1}
\end{theo}

\section{Automorphisms of the Lie algebra \boldmath{$S(M)$}}

We will finish the  description  of  automorphisms  of  the  Lie
algebra $S(M)$. We have already reduced the problem to
automorphisms which  are identity on  $\A=C^\infty(M)$.  Such an
automorphism,  respecting  the filtration, restricts to an
automorphism of $\D^1(M)=S^1(M)$, where,  in view of Theorem
\ref{d1}, it  is  of  the  form  $\Phi(X+f)=X+(f+\lambda
div(X)+\omega(X))$.       Using       the        automorphism
$e^{\overline{\omega}}$, we can reduce to the case when
$\omega=0$.  We will show that in this case $\lambda=0$ and
$\Phi=id$.

Consider an automorphism $\Phi$ of $S$ which is identical  on
functions and of the form $\Phi(X)=X+\lambda divX$ on vector
fields. It is  easy to see that this implies that
$\Phi=id_S+\psi$, where $\psi:S\raa S$  is lowering.   The
automorphism-property   yields    $\psi(\{    P,f\})
=\{\psi(P),f\},$ for all $P\in S$, $f\in\A$. Let us take $P\in
S_2$. Then $\{ P,f\}$ is a vector field (linear function on
$T^*M$) and we get
\[\lambda div\{ P,f\}=\{\psi(P),f\}.\]
But for $\lambda\ne 0,$ the l.h.s. is a second order differential
operator  with  respect  to  $f$  (e.g. for $P=X^2,$ the principal
symbol is $2\lambda X^2$), while  the r.h.s. is of first order; a
contradiction. Thus $\lambda=0$ and $\Phi$ is identity on
first-order operators (polynomials).

Now, we can proceed inductively, showing that $\psi\!\mid_{S_i}=0$
also for $i>1$. For any $P\in S,$ we have
\begin{equation}\label{ind}
\{\psi(P),f\}=\psi(\{ P,f\})\quad\textit{and}\quad
\{\psi(P),X\}=\psi(\{ P,X\}),
\end{equation}
for   any   function   $f$   and   any   vector   field    $X$.
Then, $\psi\!\mid_{S_{i-1}}=0$ and (\ref{ind}) imply that, for
$P\in S_i,$  we have $\psi(P)\in\A$ and that $\psi$ is an
intertwining operator for  the action of vector fields on $S_i$
and $\A$. But following the methods of \cite{NP} or \cite{BHMP},
one can easily see that such operators are trivial, so
$\psi\!\mid_{S_i}=0$. Thus we get $\psi=0$, i.e.  $\Phi=id_S$, and
we can formulate the following final result.
\begin{theo} A linear map $\Phi:S(M)\raa S(M)$ is an automorphism
of the Lie algebra $S(M)$ of polynomial functions on $T^*M$ with respect
to the canonical symplectic bracket if and only if it can be written
in the form
\begin{equation}\label{sm}
\Phi(P)=U_{\kappa}(P)\circ\phi^*\circ Exp(\omega^v),
\end{equation}
where $\kappa$ is a non-zero  constant,
$U_{\kappa}(P)=\kappa^{1-i}P$ for $P\in  S_i$, $\phi^*$ is the
phase lift  of  a diffeomorphism $\phi$  of  $M$  and
$Exp(\omega^v)$ is the vertical  symplectic diffeomorphism of
$T^*M$ being the translation by  a closed 1-form  $\omega$  on
$M$.  All  the objects $\kappa,\phi,\omega$ are uniquely
determined by $\Phi$. \label{SM}\end{theo} The automorphisms of
the whole  Poisson  algebra  $C^\infty(N)$ on a symplectic (or
even a Poisson)  manifold  $N$, have  been described  in
\cite{AG,JG2}. Our symplectic manifold  is particular  here (e.g.
$N=T^*M$ is non-compact and the symplectic form is exact), so  the
result of \cite{AG}  says  that automorphisms  of  the  Poisson
algebra $C^\infty(T^*M)$ are of the form $P\mapsto
sP\circ\tilde\phi$, where $s$ is   a   non-zero constant   and
$\tilde\phi$   is a conformal symplectomorphism with the conformal
constant $s$. In our case, we
deal with a subalgebra  of polynomial functions which  is
preserved  only by two   types of   conformal symplectomorphisms:
phase lifts of diffeomorphisms of $M$  and vertical
symplectomorphisms associated with closed 1-forms on $M$. In both
cases we have symplectomorphisms, so $s=1$. So far so good, the
pictures coincide,  but for $S$ we get an additional  family of
automorphisms $U_{\kappa}$. These automorphisms simply  do  not
extend  to automorphisms of the whole   algebra $C^\infty(T^*M)$.

\section{Automorphisms of the Lie algebra \boldmath{$\D(M)$}}

Let us go back to the general problem of the determination of the
automorphisms $\Phi_1$ of $\D=\D(M),$ with restriction
$\kappa\!\cdot\! id$ ($\kappa\in\R,\kappa\neq 0$) on $\A=\ci(M)$.

The restriction $\Phi_1\m_{\D^1}$ has the form given by Theorem
\ref{Phi1}. When setting \begin{equation}\Phi_2=\Phi_1\circ
e^{-\kappa^{-1}\overline{\omega}},\label{phi1}\end{equation} we
obtain---as easily verified---an automorphism of $\D$, whose
restriction to $\D^1$ is
\[\Phi_2\m_{\D^1}=\left(\kappa\pi_0+(id+\beta)\circ\pi_c\right)\m_{\D^1},\mbox{ where }\beta=\lambda div.\]
In the sequel, we shall write $\Phi$ instead of $\Phi_2$ (if no
confusion is possible).
Using (\ref{H2}), one finds that \[
[[[\Phi_cD^{i}_c,f_1]_c,f_2]_c,\ldots,f_{i-1}]_c=
\kappa^{1-i}[[[D^{i}_c,f_1]_c,f_2]_c,\ldots,f_{i-1}]_c ,\] since
$\Phi_c\!\!\mid_{\D^1_c}=id$. So
$\Phi_cD^{i}_c-\kappa^{1-i}D^{i}_c\in\D^{i-1}_c$ and
\begin{equation}\Phi\!\!\mid_{\D^{i}}=
\kappa^{1-i}\mbox{ }id+\psi_i,\label{sf}\end{equation} with
$\psi_i\in {\cal L}(\D^{i},\D^{i-1})$. Notice that $\psi_0=0,$
that
$\psi_1=\left((\kappa-1)\pi_0+\beta\circ\pi_c\right)\!\!\mid_{\D^1},$
and that $\psi_if=\kappa(1-\kappa^{-i})f$.

\begin{re} Assertion (\ref{sf}) is equivalent to say that the automorphism
$\tilde{\Phi}$ of the Poisson algebra $S$ induced by $\Phi,$ is
$U_{\kappa}$ (cf. Theorem (\ref{SM})).
\end{re}
Apply now (\ref{sf}), observe that the homomorphism-property then
reads
\begin{equation}\psi_{i+j-1}[D^{i},\Delta^{j}]=\kappa^{1-j}[\psi_iD^{i},
\Delta^{j}]+
\kappa^{1-i}[D^{i},\psi_j\Delta^{j}]+[\psi_iD^{i},\psi_j\Delta^{j}],
\label{homo-filt}\end{equation} for all
$D^{i}\in\D^{i},\Delta^{j}\in\D^{j}$, and project
(\ref{homo-filt}), written for $D_c^{i}$ and $f$ ($i\ge 2$) resp.
for $D^{i}_c$ and $\Delta_c^{j}$ ($i+j\ge 3$) on $\A$:
\begin{equation}(\psi_{i,c}D^{i}_c)(f)=\kappa^{-1}\psi_{i-1,0}[D^{i}_c,f]=
(1-\kappa^{1-i})D^{i}_c(f)+\kappa^{-1}\psi_{i-1,0}[D^{i}_c,f]_c\label{H1s}
\end{equation}
and
\begin{equation}\begin{split}\psi_{i+j-1,0}[D^{i}_c,\Delta^{j}_c]&=
\left(\left(\kappa^{1-i}id+\psi_{i,c}\right)D^{i}_c\right)(\psi_{j,0}
\Delta^{j}_c)\\
    &-\left(\left(\kappa^{1-j}id+\psi_{j,c}\right)\Delta^{j}_c\right)
    (\psi_{i,0}D^{i}_c).\end{split}\label{H3s}\end{equation}
When writing (\ref{H1s}) for $i=2$, we get
\[\psi_{2,c}D^2_c=(1-\kappa^{-1})D^2_c+\kappa^{-1}\beta[D^2_c,\cdot]_c.\]
Given that $\psi_{2,c}D^2_c\in\D^1_c$, we have
\[(1-\kappa)[D^2_c,f]_c(g)=[\beta[D^2_c,\cdot]_c,f]_c(g).\]
Since $\pi_c,$ $[.,.],$ and $\beta$ are local, the same equation
holds locally. If $D^2_c=D^1_c+D^{ij}\partial_{ij}$, an easy
computation shows that
\[(1-\kappa)[D^2_c,f]_c(g)=(1-\kappa)(D^{ij}+D^{ji})\p_if\p_jg,\]
$\beta[D^2_c,f]_c=\lambda (D^{ij}+D^{ji})\p_{ij}f+...,$ where
$...$ are terms of the first order in $f$, and
\[[\beta[D^2_c,\cdot]_c,f]_c(g)=2\lambda (D^{ij}+D^{ji})\p_if\p_jg,\]
so that
\begin{equation}1-\kappa=2\lambda\label{fe1}.\end{equation}
Equation (\ref{H3s}), written for $i=1$ and $j=2$, reads
\begin{equation}\psi_{2,0}(L_X\Delta)=X(\psi_{2,0}\Delta)-\Delta(\beta
X)-\kappa^{-1}\beta [\Delta,\beta X]_c,\label{H3s2}\end{equation}
for all $X\in Vect(M)$ and all $\Delta\in\D^2_c$.

In order to show that $\psi_{2,0}\in {\cal L}(\D^2_c,\A)$ is
local, note that it follows for instance from \cite{NP} (see
section 3) that, if $D\in\D^2_c$ vanishes on an open $U\subset M$
and if $x_0\in U$, one has $D=\sum_kL_{X_k}D_k$ ($X_k\in Vect(M)$,
$D_k\in\D^2_c$), with $X_k\m_V=D_k\m_V=0$, for some neighborhood
$V\subset U$ of $x_0$. It then suffices to combine this
decomposition of $D$ and equation (\ref{H3s2}) to find that
$(\psi_{2,0}D)(x_0)=0$.

Let $U$ be a connected, relatively compact domain of local
coordinates of $M$, in which the divergence of a vector field has
the form (\ref{classdiv}).
Recall that if $\Delta\in\D^2_{c,U}$, its representation is a
polynomial $\Delta\in\ci(U,\R^n\oplus\vee^2\R^n)$. Therefore
$\psi_{2,0}\m_U\in\cl(\ci(U,\R^n\oplus\vee^2\R^n),\ci(U))_{loc}$,
with representation $\psi(\eta;\Delta)$ ($\eta\in
(\R^n)^*,\Delta\in \R^n\oplus\vee^2\R^n$).
As easily checked, equation (\ref{H3s2}) locally reads
\begin{equation}\begin{split}(X.\psi)(\eta;\Delta)&-\la
X,\eta\ra\tau_{\zeta}\psi(\eta;\Delta)+\psi(\eta+\zeta;
X\tau_{\zeta}\Delta)\\
 &-\lambda\la X,\zeta\ra\Delta(\zeta)-\kappa^{-1}\lambda^2\la
 X,\zeta\ra\left(\Delta(\eta+2\zeta)-\Delta(\eta+\zeta)
 -\Delta(\zeta)\right)=0,\label{lfe}
\end{split}\end{equation}
where $\zeta\in(\R^n)^*$ represents once more the derivatives
acting on $X$ and where $X.\psi$ is obtained by derivation of the
coefficients of $\psi$ in the direction of $X$.

Take in (\ref{lfe}) the terms of degree $0$ in $\zeta$:
\[(X.\psi)(\eta;\Delta)=0.\]
This means that the coefficients of $\psi$ are constant.

The terms of degree $1$ lead to the equation
\[\la X,\eta\ra (\zeta\p_{\eta})\psi(\eta;\Delta)
-\psi(\eta;X(\zeta\p_{\xi})\Delta)=0,\] which, if $\rho$ denotes
the natural action of $gl(n,\R)$, may be written
\[\rho(X\otimes\zeta)\left(\psi(\eta;\Delta)\right)=0.\] Note that
$\psi(\eta;\Delta)$ is completely characterized by
$\psi(\eta;Y^k)$ ($Y\in \R^n, k\in\{1,2\}$). This last expression
is a polynomial in $\eta$ and $Y$ (remark that it's homogeneous of
degree $k$ in $Y$). It follows from the description of invariant
polynomials under the action of $gl(n,\R)$ (see \cite{HW}), that
it is a polynomial in the evaluation $\la Y,\eta\ra$. Finally,
\begin{equation}\psi(\eta;Y^k)=c_k\la
Y,\eta\ra^k,\label{poleval}\end{equation} where $c_k\in \R$.

Seeking the terms of degree $2$ in $\zeta$, we find
\begin{equation}\begin{split}\frac{1}{2}\la X,\eta\ra
(\zeta\p_{\eta})^2\psi(\eta;\Delta)&-(\zeta\p_{\eta})\psi(\eta;
X(\zeta\p_{\xi})\Delta)-\frac{1}{2}\psi(\eta;X(\zeta\p_{\xi})^2\Delta)\\
 &=-\lambda\la X,\zeta\ra\Delta^1(\zeta)-\kappa^{-1}\lambda^2\la X,\zeta\ra
 \left((\zeta\p_{\eta})\Delta(\eta)-\Delta^1(\zeta)\right),
\end{split}\label{z2}\end{equation} where $\Delta^1(\zeta)$ denotes the
terms of degree
$1$ in $\Delta(\zeta)$. Substitute now $Y^k$ ($k\in\{1,2\}$) to
$\Delta$ and observe that $X(\zeta\p_{\xi})Y^k=\la
Y,\zeta\ra(X\p_Y)Y^k$ and $X(\zeta\p_{\xi})^2Y^k=k\la Y,\zeta\ra
^2(X\p_Y)Y^{k-1}$. The l.h.s. of (\ref{z2}) then reads
\[\frac{1}{2}\la X,\eta\ra
(\zeta\p_{\eta})^2\psi(\eta;Y^k)-\la
Y,\zeta\ra(\zeta\p_{\eta})(X\p_Y)\psi(\eta;Y^k)-\frac{k}{2}\la
Y,\zeta\ra ^2(X\p_Y)\psi(\eta;Y^{k-1}).\] When setting $k=1$, then
$k=2$, when using (\ref{poleval}) (if $k=1$, the last term of the
l.h.s. vanishes) and noticing that the evaluations $\la
X,\eta\ra,$ $\la X,\zeta\ra,$ $\la Y,\eta\ra,$ and $\la
Y,\zeta\ra$ can be viewed, if $n>1$, as independent variables, one
gets from equation (\ref{z2})
\begin{equation}
c_1=\lambda,c_1+c_2=0,c_2=\kappa^{-1}\lambda^2.\label{fe2}\end{equation}
If $n=1$, one only finds $c_1=\lambda$ and
$c_1+3c_2=2\kappa^{-1}\lambda^2,$ but when selecting in
(\ref{lfe}) the terms of degree $3$ in $\zeta$, one gets
$c_1+2c_2=\lambda+2\kappa^{-1}\lambda^2,$ so that (\ref{fe2})
still holds. The solutions of the system (\ref{fe1}), (\ref{fe2})
are
$\kappa=1,\lambda=0,c_1=c_2=0$ and $\kappa=-1,\lambda=1,c_1=1,c_2=-1$.\\

Let us first examine the case $\kappa=1$. Equation (\ref{H3s}),
written---more generally---for $i=1$ and $j\ge 2$, reads
\[\psi_{j,0}(L_X\Delta^{j})=L_X(\psi_{j,0}\Delta^{j}),
\forall X\in Vect(M),\forall\Delta^{j}\in\D^{j},\] since
$\psi_{j}\!\!\mid_{\A}=0$, so that $\psi_{j,0}$ is an intertwining
operator from $(\D^{j},L)$ into $(\D^0,L)$. The results of
\cite{NP} or \cite{BHMP} show that
$\psi_{j,0}=\lambda_j\pi_0\!\!\mid_{\D^{j}}$ ($\lambda_j\in \R$),
for all $j\ge 2$ (and all $n\ge 1$; indeed, a straightforward
adaptation of the method of \cite{NP} immediately shows that this
particular result is also valid in dimension $n=1$). It's now easy
to verify that $\Phi_2$ (see (\ref{phi1})) is $id$
and that $\Phi_1=e^{\overline{\omega}}$.\\

If $\kappa=-1$, one has $\psi_{1,0}\m_{\D^1_c}=\C_0\m_{\D^1_c}$,
where $\C$ is the automorphism introduced in Lemma \ref{C}.

Inductively, if $\psi_{j-1,0}\m_{\D^{j-1}_c}=\C_0\m_{\D^{j-1}_c}$
($j\ge 2$), the same relation holds for $j$. Indeed, we obtain
from (\ref{H1s}) and (\ref{H3s}),
\begin{equation}\psi_{j,0}(L_X\Delta)=L_X(\psi_{j,0}\Delta)-\Delta(\C_0X)
+\C_0[\Delta,\C_0X]_c,\label{induc}\end{equation}
for all $X\in Vect(M)$ and all $\Delta\in\D^{j}_c$.
Straightforward computations, using the properties of $\C$, show
that
\begin{equation}\C_0[\Delta,\cdot]_c=\Delta-\C_c\Delta=\Delta-\C\Delta
+\C_0\Delta\label{C2}\end{equation}
on $\A$, as $[\Delta,f]_c=[\Delta,f]-\Delta(f)$ ($f\in\A$), and
that
\begin{equation}\C_0(L_X\Delta)=L_X(\C_0\Delta)-(\C_c\Delta)(\C_0X).
\label{C3}\end{equation}
It follows from (\ref{induc}),(\ref{C2}), and (\ref{C3}), that
\[\psi_{j,0}(L_X\Delta)-L_X(\psi_{j,0}\Delta)=\C_0(L_X\Delta)-
L_X(\C_0\Delta).\] This last equation is still valid for
$\Delta\in\D^{j}$ and signifies that $\psi_{j,0}-\C_0\m_{\D^{j}}$
is an intertwining operator from $(\D^{j},L)$ into $(\D^0,L)$.
Thus, applying once more the results of \cite{NP} or \cite{BHMP},
one sees that
\begin{equation}\psi_{j,0}\m_{\D^{j}_c}=\C_0\m_{\D^{j}_c}\label{psij0}.
\end{equation}
It is now again easy to prove (use (\ref{sf}) on $\D^{j}_c$,
(\ref{psij0}), (\ref{H1s}), and (\ref{C2})) that $\Phi_2=\C$ and
that $\Phi_1=\C\circ e^{\overline{\omega}}$.
Hence the only automorphisms of $\D$ that coincide with
$\kappa\!\cdot\! id$ ($\kappa\in \R,\kappa\neq 0$) on functions, are
$\Phi_1=e^{\overline{\omega}}$ (here $\kappa=1$) and
$\Phi_1=\C\circ e^{\overline{\omega}}$ (here $\kappa=-1$), where
$\omega$ is a closed $1$-form on $M$. Summarizing, we get the
following characterization.
\begin{theo} A linear map $\Phi:\D(M)\raa\D(M)$ is  an
automorphism of the Lie  algebra $\D(M)$ of linear differential
operators on $C^\infty(M)$ if and only if it can be written
in the form
\begin{equation} \Phi=\phi_*\circ\C^a\circ e^{\overline{\omega}},
\end{equation}
where $\phi$  is  a  diffeomorphism  of  $M$, $a=0,1$, $\C^0=id$
and $\C^1=\C$, and $\omega$  is  a  closed 1-form  on $M$. All the
objects $\phi,a,\omega$ are uniquely determined by $\Phi$.
\end{theo}
Let us notice finally that the above theorem states once more that
all automorphisms of $\D(M)$ respect the filtration and thus shows
that one-parameter groups of automorphisms of the Lie algebra
$\D(M)$ (for any reasonable topology on $\D(M)$) cannot have as
generators the inner derivations $ad_D$ for $D$ not being of the
first order. An analogous fact holds for the Lie algebra $S(M)$.
Thus we have the following.
\begin{cor}
The Lie algebras $\D(M)$ and $S(M)$ of linear differential
operators on $C^\infty(M)$ resp. of the principal symbols of these
operators, are not integrable, i.e. there are no
(infinite-dimensional) Lie groups for which they are the Lie
algebras.
\end{cor}

\noindent Janusz GRABOWSKI\\Polish Academy of Sciences, Institute
of Mathematics\\ul. 'Sniadeckich 8, P.O. Box 137\\00-950 Warsaw,
Poland\\E-mail: jagrab@impan.gov.pl\\\\
\noindent Norbert PONCIN\\Universit\'e de Luxembourg, D\'epartement de Math\'ematiques\\Avenue de la Fa\"{\i}encerie, 162 A\\
L-1511 Luxembourg, Luxembourg\\E-mail: poncin@cu.lu

\end{document}